\documentclass[reqno,12pt,a4paper]{amsart}

\usepackage{mathrsfs}
\usepackage{amssymb}
\usepackage{amsfonts}
\usepackage{latexsym}
\usepackage{amsthm}
\usepackage{graphicx}

\def\lb{\label}

\newcommand{\er}[1]{\textrm{(\ref{#1})}}

\begin{document}


\renewcommand{\theequation}{\arabic{section}.\arabic{equation}}
\theoremstyle{plain}
\newtheorem{theorem}{\bf Theorem}[section]
\newtheorem{lemma}[theorem]{\bf Lemma}
\newtheorem{corollary}[theorem]{\bf Corollary}
\newtheorem{proposition}[theorem]{\bf Proposition}
\newtheorem{definition}[theorem]{\bf Definition}
\newtheorem{example}[theorem]{\bf Example}

\newtheorem{remark}[theorem]{\bf Remark}

\def\a{\alpha}  \def\cA{{\mathcal A}}     \def\bA{{\bf A}}  \def\mA{{\mathscr A}}
\def\b{\beta}   \def\cB{{\mathcal B}}     \def\bB{{\bf B}}  \def\mB{{\mathscr B}}
\def\g{\gamma}  \def\cC{{\mathcal C}}     \def\bC{{\bf C}}  \def\mC{{\mathscr C}}
\def\G{\Gamma}  \def\cD{{\mathcal D}}     \def\bD{{\bf D}}  \def\mD{{\mathscr D}}
\def\d{\delta}  \def\cE{{\mathcal E}}     \def\bE{{\bf E}}  \def\mE{{\mathscr E}}
\def\D{\Delta}  \def\cF{{\mathcal F}}     \def\bF{{\bf F}}  \def\mF{{\mathscr F}}
\def\c{\chi}    \def\cG{{\mathcal G}}     \def\bG{{\bf G}}  \def\mG{{\mathscr G}}
\def\z{\zeta}   \def\cH{{\mathcal H}}     \def\bH{{\bf H}}  \def\mH{{\mathscr H}}
\def\e{\eta}    \def\cI{{\mathcal I}}     \def\bI{{\bf I}}  \def\mI{{\mathscr I}}
\def\p{\psi}    \def\cJ{{\mathcal J}}     \def\bJ{{\bf J}}  \def\mJ{{\mathscr J}}
\def\vT{\Theta} \def\cK{{\mathcal K}}     \def\bK{{\bf K}}  \def\mK{{\mathscr K}}
\def\k{\kappa}  \def\cL{{\mathcal L}}     \def\bL{{\bf L}}  \def\mL{{\mathscr L}}
\def\l{\lambda} \def\cM{{\mathcal M}}     \def\bM{{\bf M}}  \def\mM{{\mathscr M}}
\def\L{\Lambda} \def\cN{{\mathcal N}}     \def\bN{{\bf N}}  \def\mN{{\mathscr N}}
\def\m{\mu}     \def\cO{{\mathcal O}}     \def\bO{{\bf O}}  \def\mO{{\mathscr O}}
\def\n{\nu}     \def\cP{{\mathcal P}}     \def\bP{{\bf P}}  \def\mP{{\mathscr P}}
\def\r{\rho}    \def\cQ{{\mathcal Q}}     \def\bQ{{\bf Q}}  \def\mQ{{\mathscr Q}}
\def\s{\sigma}  \def\cR{{\mathcal R}}     \def\bR{{\bf R}}  \def\mR{{\mathscr R}}
\def\S{\Sigma}  \def\cS{{\mathcal S}}     \def\bS{{\bf S}}  \def\mS{{\mathscr S}}
\def\t{\tau}    \def\cT{{\mathcal T}}     \def\bT{{\bf T}}  \def\mT{{\mathscr T}}
\def\f{\phi}    \def\cU{{\mathcal U}}     \def\bU{{\bf U}}  \def\mU{{\mathscr U}}
\def\F{\Phi}    \def\cV{{\mathcal V}}     \def\bV{{\bf V}}  \def\mV{{\mathscr V}}
\def\P{\Psi}    \def\cW{{\mathcal W}}     \def\bW{{\bf W}}  \def\mW{{\mathscr W}}
\def\o{\omega}  \def\cX{{\mathcal X}}     \def\bX{{\bf X}}  \def\mX{{\mathscr X}}
\def\x{\xi}     \def\cY{{\mathcal Y}}     \def\bY{{\bf Y}}  \def\mY{{\mathscr Y}}
\def\X{\Xi}     \def\cZ{{\mathcal Z}}     \def\bZ{{\bf Z}}  \def\mZ{{\mathscr Z}}
\def\be{{\bf e}}
\def\bv{{\bf v}} \def\bu{{\bf u}}
\def\Om{\Omega}
\def\vr{\varrho}
\def\bbD{\pmb \Delta}
\def\mm{\mathrm m}
\def\mn{\mathrm n}

\newcommand{\mc}{\mathscr {c}}

\newcommand{\gA}{\mathfrak{A}}          \newcommand{\ga}{\mathfrak{a}}
\newcommand{\gB}{\mathfrak{B}}          \newcommand{\gb}{\mathfrak{b}}
\newcommand{\gC}{\mathfrak{C}}          \newcommand{\gc}{\mathfrak{c}}
\newcommand{\gD}{\mathfrak{D}}          \newcommand{\gd}{\mathfrak{d}}
\newcommand{\gE}{\mathfrak{E}}
\newcommand{\gF}{\mathfrak{F}}           \newcommand{\gf}{\mathfrak{f}}
\newcommand{\gG}{\mathfrak{G}}           
\newcommand{\gH}{\mathfrak{H}}           \newcommand{\gh}{\mathfrak{h}}
\newcommand{\gI}{\mathfrak{I}}           \newcommand{\gi}{\mathfrak{i}}
\newcommand{\gJ}{\mathfrak{J}}           \newcommand{\gj}{\mathfrak{j}}
\newcommand{\gK}{\mathfrak{K}}            \newcommand{\gk}{\mathfrak{k}}
\newcommand{\gL}{\mathfrak{L}}            \newcommand{\gl}{\mathfrak{l}}
\newcommand{\gM}{\mathfrak{M}}            \newcommand{\gm}{\mathfrak{m}}
\newcommand{\gN}{\mathfrak{N}}            \newcommand{\gn}{\mathfrak{n}}
\newcommand{\gO}{\mathfrak{O}}
\newcommand{\gP}{\mathfrak{P}}             \newcommand{\gp}{\mathfrak{p}}
\newcommand{\gQ}{\mathfrak{Q}}             \newcommand{\gq}{\mathfrak{q}}
\newcommand{\gR}{\mathfrak{R}}             \newcommand{\gr}{\mathfrak{r}}
\newcommand{\gS}{\mathfrak{S}}              \newcommand{\gs}{\mathfrak{s}}
\newcommand{\gT}{\mathfrak{T}}             \newcommand{\gt}{\mathfrak{t}}
\newcommand{\gU}{\mathfrak{U}}             \newcommand{\gu}{\mathfrak{u}}
\newcommand{\gV}{\mathfrak{V}}             \newcommand{\gv}{\mathfrak{v}}
\newcommand{\gW}{\mathfrak{W}}             \newcommand{\gw}{\mathfrak{w}}
\newcommand{\gX}{\mathfrak{X}}               \newcommand{\gx}{\mathfrak{x}}
\newcommand{\gY}{\mathfrak{Y}}              \newcommand{\gy}{\mathfrak{y}}
\newcommand{\gZ}{\mathfrak{Z}}             \newcommand{\gz}{\mathfrak{z}}

\def\ve{\varepsilon}   \def\vt{\vartheta}    \def\vp{\varphi}    \def\vk{\varkappa}

\def\A{{\mathbb A}} \def\B{{\mathbb B}} \def\C{{\mathbb C}}
\def\dD{{\mathbb D}} \def\E{{\mathbb E}} \def\dF{{\mathbb F}} \def\dG{{\mathbb G}} \def\H{{\mathbb H}}\def\I{{\mathbb I}} \def\J{{\mathbb J}} \def\K{{\mathbb K}} \def\dL{{\mathbb L}}\def\M{{\mathbb M}} \def\N{{\mathbb N}} \def\O{{\mathbb O}} \def\dP{{\mathbb P}} \def\R{{\mathbb R}}\def\S{{\mathbb S}} \def\T{{\mathbb T}} \def\U{{\mathbb U}} \def\V{{\mathbb V}}\def\W{{\mathbb W}} \def\X{{\mathbb X}} \def\Y{{\mathbb Y}} \def\Z{{\mathbb Z}}


\def\la{\leftarrow}              \def\ra{\rightarrow}            \def\Ra{\Rightarrow}
\def\ua{\uparrow}                \def\da{\downarrow}
\def\lra{\leftrightarrow}        \def\Lra{\Leftrightarrow}


\def\lt{\biggl}                  \def\rt{\biggr}
\def\ol{\overline}               \def\wt{\widetilde}
\def\ul{\underline}
\def\no{\noindent}


\let\ge\geqslant                 \let\le\leqslant
\def\lan{\langle}                \def\ran{\rangle}
\def\/{\over}                    \def\iy{\infty}
\def\sm{\setminus}               \def\es{\emptyset}
\def\ss{\subset}                 \def\ts{\times}
\def\pa{\partial}                \def\os{\oplus}
\def\om{\ominus}                 \def\ev{\equiv}
\def\iint{\int\!\!\!\int}        \def\iintt{\mathop{\int\!\!\int\!\!\dots\!\!\int}\limits}
\def\el2{\ell^{\,2}}             \def\1{1\!\!1}
\def\sh{\sharp}
\def\wh{\widehat}
\def\bs{\backslash}
\def\intl{\int\limits}

\def\na{\mathop{\mathrm{\nabla}}\nolimits}
\def\sh{\mathop{\mathrm{sh}}\nolimits}
\def\ch{\mathop{\mathrm{ch}}\nolimits}
\def\where{\mathop{\mathrm{where}}\nolimits}
\def\all{\mathop{\mathrm{all}}\nolimits}
\def\as{\mathop{\mathrm{as}}\nolimits}
\def\Area{\mathop{\mathrm{Area}}\nolimits}
\def\arg{\mathop{\mathrm{arg}}\nolimits}
\def\const{\mathop{\mathrm{const}}\nolimits}
\def\det{\mathop{\mathrm{det}}\nolimits}
\def\diag{\mathop{\mathrm{diag}}\nolimits}
\def\diam{\mathop{\mathrm{diam}}\nolimits}
\def\dim{\mathop{\mathrm{dim}}\nolimits}
\def\dist{\mathop{\mathrm{dist}}\nolimits}
\def\Im{\mathop{\mathrm{Im}}\nolimits}
\def\Iso{\mathop{\mathrm{Iso}}\nolimits}
\def\Ker{\mathop{\mathrm{Ker}}\nolimits}
\def\Lip{\mathop{\mathrm{Lip}}\nolimits}
\def\rank{\mathop{\mathrm{rank}}\limits}
\def\Ran{\mathop{\mathrm{Ran}}\nolimits}
\def\Re{\mathop{\mathrm{Re}}\nolimits}
\def\Res{\mathop{\mathrm{Res}}\nolimits}
\def\res{\mathop{\mathrm{res}}\limits}
\def\sign{\mathop{\mathrm{sign}}\nolimits}
\def\span{\mathop{\mathrm{span}}\nolimits}
\def\supp{\mathop{\mathrm{supp}}\nolimits}
\def\Tr{\mathop{\mathrm{Tr}}\nolimits}
\def\BBox{\hspace{1mm}\vrule height6pt width5.5pt depth0pt \hspace{6pt}}


\newcommand\nh[2]{\widehat{#1}\vphantom{#1}^{(#2)}}
\def\dia{\diamond}

\def\Oplus{\bigoplus\nolimits}



\def\qqq{\qquad}
\def\qq{\quad}
\let\ge\geqslant
\let\le\leqslant
\let\geq\geqslant
\let\leq\leqslant
\newcommand{\ca}{\begin{cases}}
\newcommand{\ac}{\end{cases}}
\newcommand{\ma}{\begin{pmatrix}}
\newcommand{\am}{\end{pmatrix}}
\renewcommand{\[}{\begin{equation}}
\renewcommand{\]}{\end{equation}}
\def\eq{\begin{equation}}
\def\qe{\end{equation}}
\def\[{\begin{equation}}
\def\bu{\bullet}

\title[{Schr\"odinger operators with guided potentials on periodic graphs}]
{Schr\"odinger operators  with guided potentials on periodic graphs}

\date{\today}
\author[Evgeny Korotyaev]{Evgeny Korotyaev}
\address{Saint-Petersburg State University, Universitetskaya nab. 7/9, St. Petersburg, 199034, Russia,
\ korotyaev@gmail.com, \
e.korotyaev@spbu.ru,}
\author[Natalia Saburova]{Natalia Saburova}
\address{Northern (Arctic) Federal University, Severnaya Dvina Emb. 17,
Arkhangelsk, 163002, Russia,
 \ n.saburova@gmail.com, \ n.saburova@narfu.ru}


\begin{abstract}

We consider discrete Schr\"odinger operators with periodic potentials on
periodic graphs perturbed by guided non-positive potentials, which
are periodic in some directions and finitely supported in other
ones. The spectrum of the unperturbed operator is a union of a
finite number of non-degenerate bands and eigenvalues of infinite
multiplicity. We show that the spectrum of the perturbed operator
consists of the "unperturbed" one plus the additional guided
spectrum, which is a union of a finite number of bands. We estimate
the position of the guided bands and their length in terms of graph
geometric parameters. We also determine the asymptotics of the guided
bands for large guided potentials. Moreover, we show that the possible number of the guided bands, their length and position can be rather arbitrary for
some specific potentials.
\end{abstract}

\maketitle

\vskip 0.25cm

\section {\lb{Sec1}Introduction}
\setcounter{equation}{0} Discrete Schr\"odinger
operators on periodic graphs have attracted a lot of attention due
to their applications to the study of electronic properties of real
crystalline structures, see, e.g., \cite{Ha02},
\cite{H89}, \cite{NG04} and the survey \cite{CGPNG09}.
Waveguides defects  allow one to obtain conductivity of the
material for those frequencies (energies) at which it was not in
purely periodic structure.  Such effects have a lot of applications,
see about waveguides in photonic crystal structures in \cite{J00},
\cite{JJ02}, \cite{H15} and references therein.

We consider discrete Schr\"odinger operators with periodic potentials on
periodic graphs perturbed by guided non-positive potentials, which
are periodic in some directions and finitely supported in other
ones. For example, on the lattice $\Z^2$ the support of the guided
potentials is a strip. It is well-known that the spectrum of
Schr\"odinger operators with periodic potentials on periodic graphs
has a band structure with a finite number of flat bands (eigenvalues
of infinite multiplicity) \cite{HN09}, \cite{HS04}, \cite{KS14}, \cite{RR07}. The
spectrum of the perturbed Schr\"odinger operator consists of the
spectrum of the "unperturbed" operator plus the \emph{guided}
spectrum. The additional guided spectrum is a union of a finite
number of bands, here the corresponding wave-functions are located
along the support of the guided potentials and decrease in
perpendicular directions. Note that line defects on the lattice were
considered in \cite{C12}, \cite{Ku14}, \cite{Ku16}, \cite{OA12}.

\medskip

In our paper we study the  influence of  guided potentials on the
spectrum of Schr\"odinger operators. We describe our main goals:

1) to estimate the position and the length of the guided bands in terms of geometric
parameters of graphs and guided potentials;

2) to determine the asymptotics of the guided spectrum  for large guided
potentials;

3) to show that the possible number of the guided bands, their length and position can be rather arbitrary for
some specific potentials.

\subsection{Schr\"odinger operators with periodic potentials}
Let $\G=(V,\cE)$ be a connected infinite graph, possibly  having
loops and multiple edges, where $V$ is the set of its vertices and
$\cE$ is the set of its unoriented edges. From the set $\cE$ we
construct the set $\cA$ of oriented edges by considering each edge
in $\cE$ to have two orientations. An edge starting at a vertex $u$
and ending at a vertex $v$ from $V$ will be denoted as the ordered
pair $(u,v)\in\cA$.
We define the degree
${\vk}_v$ of the vertex $v\in V$ as the number of all edges from
$\cA$ starting at $v$.

Below we consider locally finite
$\Z^{\wt d}$-periodic graphs $\G$, $\wt d\geq2$, i.e., graphs satisfying the
following conditions:
\begin{itemize}
  \item[1)] \emph{$\G$ is equipped with an action of the free abelian group $\Z^{\wt d}$;}
  \item[2)] \emph{the degree of each vertex is finite;}

  \item[3)] \emph{the quotient graph  $\G_*=(V_*,\cE_*)=\G/\Z^{\wt d}$ is finite.}
\end{itemize}

We assume that the graphs are embedded into Euclidean space, since
in many applications such a natural embedding exists. For example,
in the tight-binding approximation real crystalline structures are
modeled as discrete graphs embedded into $\R^d$ ($d=2,3$) and
consisting of vertices (points representing positions of atoms) and
edges (representing chemical bonding of atoms), by ignoring the
physical characters of atoms and bonds that may be different from
one another. But all results of the paper stay valid in the case
of abstract periodic graphs (without the assumption of graph embedding
into Euclidean space).

For a periodic graph $\G$ embedded into the space $\R^{\wt d}$, the
quotient graph  $\G/\Z^{\wt d}$ is a graph on the $\wt
d$-dimensional torus $\R^{\wt d}/\Z^{\wt d}$. Due to the definition, the graph $\G$ is invariant under translations through vectors $a_1,\ldots,a_{\wt d}$\, which generate the group $\Z^{\wt d}$:
$$
\G+a_s=\G,\qqq \forall\, s\in\N_{\wt d}\,.
$$
Here and below for each integer
$m$ the set $\N_m$ is given by
\[
\N_{m}=\{1,\ldots,m\,\}.
\]
We will call the vectors $a_1,\ldots,a_{\wt d}$ \emph{the periods of the graph} $\G$. In the space $\R^{\wt
d}$ we consider a coordinate system with the origin at some point
$O$ and with the basis $a_1,\ldots,a_{\wt d}$. Below the coordinates
of all vertices of $\G$ will be expressed  in this coordinate
system.

\

Let $\ell^2(V)$ be the
Hilbert space of all functions $f:V\to \C$ equipped
with the norm
$$
\|f\|^2_{\ell^2(V)}=\sum_{v\in V}|f(v)|^2<\infty.
$$
For a self-adjoint operator $A$, $\sigma(A)$,  $\sigma_{ess}(A)$, $\sigma_{ac}(A)$,
$\sigma_{p}(A)$, and $\sigma_{fb}(A)$ denote its
spectrum,  essential spectrum, absolutely continuous spectrum, point spectrum (eigenvalues of finite multiplicity), and the set of all its flat bands (eigenvalues of infinite multiplicity), respectively.

We consider  a discrete
Schr\"odinger operator $H_0$ with a periodic potential $W$ on
$f\in\ell^2(V)$ as \emph{an unperturbed operator} defined  by
\[ \lb{Sh} H_0=\D+W,
\]
where $\D$ is the discrete Laplacian (i.e., the combinatorial
 Laplace operator) given by
\[
\lb{DOLN}
\begin{aligned}
\big(\D f\big)(v)=\sum_{(v,u)\in\cA}\big(f(v)-f(u)\big),
 \qquad \forall\,f\in\ell^2(V), \qqq \forall\, v\in V,
 \end{aligned}
\]
and the sum in \er{DOLN} is taken over all oriented edges starting
at the  vertex $v\in V$. The potential $W$ is real valued and
$\Z^{\wt d}$-periodic, i.e.,
\[
\lb{potW}
(Wf)(v)=W(v)f(v),\qqq W(v+a_s)=W(v), \qqq  \forall\,(v,s)\in V\ts\N_{\wt d}\,.
\]
It is well-known that $H_0$ is self-adjoint and its
spectrum is a union of $\n$ spectral bands $\s_n(H_0)$:
\[
\lb{sH0}
\s(H_0)=\bigcup_{n=1}^\n\s_n(H_0)=\s_{ac}(H_0)\cup\s_{fb}(H_0),
\]
where $\n=\# V_*$ is the number of vertices of the quotient graph
$\G_*$, the absolutely continuous spectrum $\s_{ac}(H_0)$ consists
of non-degenerate bands $\s_n(H_0)$; $\s_{fb}(H_0)$ is the set of
all flat bands. Without loss of generality assume that the spectrum
$\s(H_0)$ is a subset of the interval $[0,\vr]$:
\[
\lb{mm} \s(H_0)\ss [0,\vr],\qqq \inf \s(H_0)=0,\qqq \vr=\sup
\s(H_0).
\]

\subsection{Results overview}
There are results about spectral properties of the discrete
Schr\"odinger operator $H_0$ with a periodic potential  $W$. The
decomposition of the operator $H_0$ into a constant fiber direct
integral was obtained in \cite{HN09}, \cite{HS04}, \cite{RR07}
without an exact form of fiber operators and in \cite{KS14},
\cite{KS17} with an exact form of fiber operators. In particular,
this yields the band-gap structure of the spectrum of the operator $H_0$. In \cite{GKT93} the authors considered the
 Schr\"odinger operator with a periodic potential on the lattice $\Z^2$, the
simplest example of $\Z^2$-periodic graphs.  They studied its Bloch
variety and its integrated density of states. In \cite{LP08},
\cite{KS15} the positions of the spectral bands of the Laplacians
were estimated in terms of eigenvalues of the operator on finite
graphs (the so-called eigenvalue bracketing).
The estimate of the total length of all spectral bands
$\s_n(H_0)$
\[
\lb{eq.7}
\sum_{n=1}^{\n}|\s_n(H_0)|\le 2\b,
\]
was obtained in \cite{KS14}; where $\b=\#\cE_*-\n+1$ is the
so-called Betti number, $\#\cE_*$ is the number of edges of the
quotient graph $\G_*$. Moreover,
 a global variation of the Lebesgue measure of the spectrum
and a global variation of the gap-length in terms of potentials and
geometric parameters of the graph were determined. Note that the
estimate \er{eq.7} also holds true for magnetic Schr\"odinger
operators with periodic magnetic and electric potentials (see
\cite{KS17}). Estimates of the Lebesgue measure of the spectrum of
$H_0$ in terms of eigenvalues of Dirichlet and Neumann operators on
a fundamental domain of the periodic graph were described in
\cite{KS15}. Estimates of effective masses, associated with the
ends of each spectral band of the Laplacian, in terms of geometric parameters of the
graphs were obtained in \cite{KS16}. Moreover, in the case of
the bottom of the spectrum two-sided estimates on the effective mass
in terms of geometric parameters of the graphs were determined. The
proof of all these results in \cite{KS14}-\cite{KS17} is based on
Floquet theory and the exact form of fiber Schr\"odinger operators
from \cite{KS14}.  The spectra of the discrete Schr\"odinger
operators on graphene nano-tubes and nano-ribbons in external fields
were discussed in \cite{KK10}, \cite{KK10a}. Finally, we note that
different properties of Schr\"odinger operators on graphs  were
considered in \cite{G15}, \cite{Sh98}.

Scattering theory for self-adjoint Schr\"odinger operators with
decreasing potentials  was investigated in \cite{BS99}, \cite{IK12}
(for the lattice) and in \cite{PR16} (for  periodic graphs). Inverse
scattering theory with finitely supported potentials was considered
in \cite{IK12} for the case of the lattice $\Z^d$ and in \cite{A12}
for the case of the hexagonal lattice. The absence of eigenvalues
embedded in the essential spectrum of the operators was discussed in
\cite{IM14}, \cite{V14}. Trace formulae and global eigenvalues
estimates for Schr\"odinger operators with complex decaying
potentials on the lattice were obtained in \cite{KL16}.
The Cwikel-Lieb-Rosenblum type bound for the discrete Schr{\"o}dinger
operator on $\Z^d$ was computed in \cite{Ka08}, \cite{RS09}.

\section {\lb{Sec1.1}Main results}
\setcounter{equation}{0}
\subsection{Schr\"odinger operators with guided potentials}
Let integer $d<\wt d$. We define the infinite \emph{fundamental graph} $\cC=\G/\Z^d$ of the
$\Z^{\wt d}$-periodic graph $\G$, which is a graph on the cylinder $\R^{\wt
d}/\Z^d$.  We also call the fundamental graph $\cC$ \emph{a discrete
cylinder} or just \emph{a cylinder}. The cylinder
$\cC=(V_c,\cE_c)$ has the vertex set $V_c=V/\Z^d$, the set
$\cE_c=\cE/\Z^d$ of unoriented edges and the set $\cA_c=\cA/\Z^d$ of
oriented edges. Note that the fundamental graph $\cC$ is $\Z^{\wt
d-d}$-periodic.

We identify the vertices of the cylinder $\cC$ with the vertices  of
the periodic graph $\G$ from the strip $\cS=[0,1)^d\ts\R^{\wt d-d}$. We
will call this infinite vertex set \emph{the fundamental vertex set of
$\G$} and denote it by the same symbol $V_c$:
\[\lb{fvs}
V_c=V\cap\cS,\qqq \cS=[0,1)^d\ts\R^{\wt d-d}.
\]
Edges of the periodic graph $\G$ connecting the vertices from the fundamental vertex set $V_c$ with the
vertices from $V\sm V_c$ will be called \emph{bridges}. Bridges always exist and provide the connectivity of the periodic graph. The
set of all bridges of the graph $\G$ we denote by
$\cB$.

We consider  a \emph{guided Schr\"odinger operator} $H$ on the
periodic graph $\G$ given by
\[\lb{gSo}
H=H_0-Q, \qqq \big(Qf\big)(v)=Q(v)f(v), \qqq f\in\ell^2(V),
\]
where $Q$ is  the \textbf{guided potential} defined by

1) {\it $Q\ge 0$ and $Q$ is   $\Z^d$-periodic, i.e.,
\[
\lb{PotP}
Q(v+a_s)=Q(v), \qqq  \forall\,(v,s)\in V_c\ts\N_d,
\]

2) the restriction of $Q$ to $V_c$ has a finite support}:
\[
\supp(Q\upharpoonright V_c)=\{v_1,\ldots,v_p\}\ss V_c\,.
\]

In other words, the guided potential $Q$ is periodic in the
directions  $a_1,\ldots,a_d$ and finitely supported in other ones,
see Fig.1.

\setlength{\unitlength}{1.0mm}
\begin{figure}[h]
\centering

\unitlength 0.7mm 
\linethickness{0.4pt}
\ifx\plotpoint\undefined\newsavebox{\plotpoint}\fi 
\begin{picture}(200,60)(0,0)

\put(-12,5){(\emph{a})}
\put(30,52){$\dL^2$}
\put(0,10){\line(1,0){85.00}}
\put(10,20){\line(1,0){85.00}}
\put(20,30){\line(1,0){85.00}}
\put(30,40){\line(1,0){85.00}}
\put(40,50){\line(1,0){85.00}}

\put(0,5){\line(1,1){50.00}}
\put(15,5){\line(1,1){50.00}}
\put(30,5){\line(1,1){50.00}}

\bezier{60}(31,5)(56,30)(81,55)
\bezier{60}(32,5)(57,30)(82,55)
\bezier{60}(33,5)(58,30)(83,55)
\bezier{60}(34,5)(59,30)(84,55)
\bezier{60}(35,5)(60,30)(85,55)
\bezier{60}(36,5)(61,30)(86,55)
\bezier{60}(37,5)(62,30)(87,55)
\bezier{60}(38,5)(63,30)(88,55)
\bezier{60}(39,5)(64,30)(89,55)
\bezier{60}(40,5)(65,30)(90,55)
\bezier{60}(41,5)(66,30)(91,55)
\bezier{60}(42,5)(67,30)(92,55)
\bezier{60}(43,5)(68,30)(93,55)
\bezier{60}(44,5)(69,30)(94,55)
\bezier{60}(45,5)(70,30)(95,55)
\bezier{60}(46,5)(71,30)(96,55)
\bezier{60}(47,5)(72,30)(97,55)
\bezier{60}(48,5)(73,30)(98,55)
\bezier{60}(49,5)(74,30)(99,55)
\bezier{60}(50,5)(75,30)(100,55)
\bezier{60}(51,5)(76,30)(101,55)
\bezier{60}(52,5)(77,30)(102,55)
\bezier{60}(53,5)(78,30)(103,55)
\bezier{60}(54,5)(79,30)(104,55)
\bezier{60}(55,5)(80,30)(105,55)
\bezier{60}(56,5)(81,30)(106,55)
\bezier{60}(57,5)(82,30)(107,55)
\bezier{60}(58,5)(83,30)(108,55)
\bezier{60}(59,5)(84,30)(109,55)
\bezier{60}(60,5)(85,30)(110,55)

\put(45,5){\line(1,1){50.00}}
\put(60,5){\line(1,1){50.00}}
\put(75,5){\line(1,1){50.00}}

\put(5,10){\circle{1}}
\put(20,10){\circle{1}}
\put(35,10){\circle*{2}}
\put(50,10){\circle*{2}}
\put(65,10){\circle{1}}
\put(80,10){\circle{1}}

\put(15,20){\circle{1}}
\put(30,20){\circle{1}}
\put(45,20){\circle*{2}}
\put(45,20){\vector(1,0){30.00}}
\put(45,20){\vector(1,1){10.00}}
\put(44,14.5){$O$}
\put(66.0,16.5){$a_1$}
\put(45.9,26.5){$a_2$}
\put(60,20){\circle*{2}}
\put(75,20){\circle{1}}
\put(90,20){\circle{1}}

\put(25,30){\circle{1}}
\put(40,30){\circle{1}}
\put(55,30){\circle*{2}}
\put(70,30){\circle*{2}}
\put(85,30){\circle{1}}
\put(100,30){\circle{1}}

\put(35,40){\circle{1}}
\put(50,40){\circle{1}}
\put(65,40){\circle*{2}}
\put(80,40){\circle*{2}}
\put(95,40){\circle{1}}
\put(110,40){\circle{1}}

\put(45,50){\circle{1}}
\put(60,50){\circle{1}}
\put(75,50){\circle*{2}}
\put(90,50){\circle*{2}}
\put(105,50){\circle{1}}
\put(120,50){\circle{1}}

\linethickness{1.5pt}
\put(15,20){\line(0,1){6.00}}
\put(30,20){\line(0,1){3.00}}
\put(45,20){\line(0,1){6.00}}
\put(60,20){\line(0,1){3.00}}
\put(75,20){\line(0,1){6.00}}
\put(90,20){\line(0,1){3.00}}

\put(25,30){\line(0,1){4.00}}
\put(40,30){\line(0,1){6.00}}
\put(55,30){\line(0,1){4.00}}
\put(70,30){\line(0,1){6.00}}
\put(85,30){\line(0,1){4.00}}
\put(100,30){\line(0,1){6.00}}
\put(35,40){\line(0,1){6.00}}
\put(50,40){\line(0,1){3.00}}
\put(65,40){\line(0,1){6.00}}
\put(80,40){\line(0,1){3.00}}
\put(95,40){\line(0,1){6.00}}
\put(110,40){\line(0,1){3.00}}
\linethickness{0.4pt}
\put(82,53){$\cS$}
\put(150,52){$\cC$}
\put(125,10){\line(1,0){30.00}}
\put(135,20){\line(1,0){30.00}}
\put(145,30){\line(1,0){30.00}}
\put(155,40){\line(1,0){30.00}}
\put(165,50){\line(1,0){30.00}}

\put(120,5){\line(1,1){50.0}}
\put(135,5){\line(1,1){50.0}}
\bezier{60}(150,5)(175,30)(200,55)
\bezier{60}(149,5)(174,30)(199,55)
\bezier{60}(148,5)(173,30)(198,55)
\bezier{60}(147,5)(172,30)(197,55)
\bezier{60}(146,5)(171,30)(196,55)
\bezier{60}(145,5)(170,30)(195,55)
\bezier{60}(144,5)(169,30)(194,55)
\bezier{60}(143,5)(168,30)(193,55)
\bezier{60}(142,5)(167,30)(192,55)
\bezier{60}(141,5)(166,30)(191,55)
\bezier{60}(140,5)(165,30)(190,55)
\bezier{60}(139,5)(164,30)(189,55)
\bezier{60}(138,5)(163,30)(188,55)
\bezier{60}(137,5)(162,30)(187,55)
\bezier{60}(136,5)(161,30)(186,55)
\bezier{60}(135,5)(160,30)(185,55)
\bezier{60}(134,5)(159,30)(184,55)
\bezier{60}(133,5)(158,30)(183,55)
\bezier{60}(132,5)(157,30)(182,55)
\bezier{60}(131,5)(156,30)(181,55)
\bezier{60}(130,5)(155,30)(180,55)
\bezier{60}(129,5)(154,30)(179,55)
\bezier{60}(128,5)(153,30)(178,55)
\bezier{60}(127,5)(152,30)(177,55)
\bezier{60}(126,5)(151,30)(176,55)
\bezier{60}(125,5)(150,30)(175,55)
\bezier{60}(124,5)(149,30)(174,55)
\bezier{60}(123,5)(148,30)(173,55)
\bezier{60}(122,5)(147,30)(172,55)
\bezier{60}(121,5)(146,30)(171,55)
\bezier{60}(120,5)(145,30)(170,55)

\put(125,10){\circle*{1.5}}
\put(140,10){\circle*{1.5}}
\put(155,10){\circle{1.5}}

\put(135,20){\circle*{1.5}}
\put(150,20){\circle*{1.5}}
\put(165,20){\circle{1.5}}

\put(145,30){\circle*{1.5}}
\put(160,30){\circle*{1.5}}
\put(175,30){\circle{1.5}}

\put(155,40){\circle*{1.5}}
\put(170,40){\circle*{1.5}}
\put(185,40){\circle{1.5}}

\put(165,50){\circle*{1.5}}
\put(180,50){\circle*{1.5}}
\put(195,50){\circle{1.5}}

\linethickness{1.5pt}
\put(135,20){\line(0,1){6.00}}
\put(150,20){\line(0,1){3.00}}
\put(145,30){\line(0,1){4.00}}
\put(160,30){\line(0,1){6.00}}
\put(155,40){\line(0,1){6.00}}
\put(170,40){\line(0,1){3.00}}
\linethickness{0.4pt}
\put(107,5){(\emph{b})}
\end{picture}

\vspace{-0.3cm} \caption{\footnotesize  \emph{a}) The square lattice $\dL^2$;
the vertices from the set $V_c$ are big black points; the strip $\cS$ is shaded;\;
\emph{b}) the discrete cylinder
$\cC=\dL^2/\Z$ (the edges of the strip are identified). The
values of the guided potential $Q$ are along the vertical axis.} \label{SLpQ}
\end{figure}
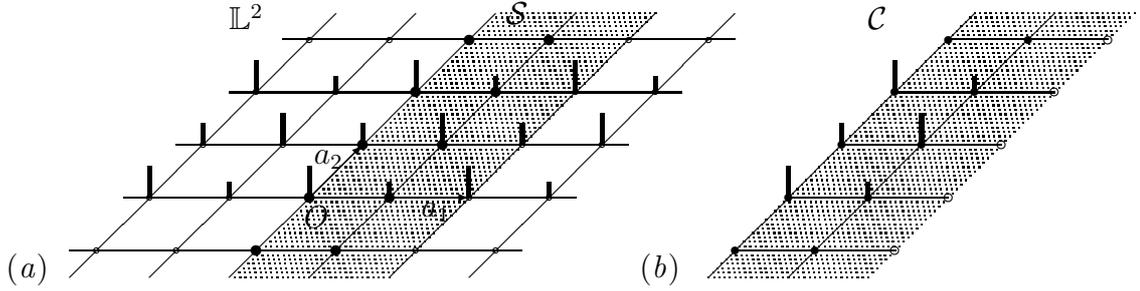

\textbf{Example.} For the square lattice $\dL^2$ with the periods $a_1,a_2$, see
Fig.\ref{SLpQ}.\emph{a}, the discrete cylinder
$\cC=\dL^2/\Z=(V_c,\cE_c)$ is
shown in Fig.\ref{SLpQ}.\emph{b}. The vertices from the set $V_c$
are big black points. The guided
potential $Q$ is shown by vertical lines. The spectrum of
a guided operator on the square lattice $\dL^2$ is discussed in Example \ref{Ex}.

We define the torus $\T^d=\R^d/(2\pi\Z)^d$ and describe the basic spectral properties of guided Schr\"odinger operators.

\begin{proposition}\label{TSDI}
i) The guided Schr\"odinger operator $H=H_0-Q$ has the following decomposition
into a constant fiber direct integral for some unitary operator $U:
\ell^2(V)\to \mH$:
\[
\lb{raz}
\begin{aligned}
& \mH=\int^\oplus_{\T^d}\ell^2(V_c){d\vt\/(2\pi)^d}\,,\qqq UH
U^{-1}=\int^\oplus_{\T^d}H(\vt){d\vt\/(2\pi)^d}\,,\qqq
H(\vt)=H_0(\vt)-Q,
\end{aligned}
\]
where the fiber Schr\"odinger operator $H(\vt)$ acts on the fiber space
$\ell^2(V_c)$ and $H_0(\vt)=\D(\vt)+W$ is the fiber operator for $H_0$, the fiber Laplacian $\D(\vt)$ is given by
\[
\label{l2.15''}
 \big(\D(\vt)f\big)(v)=\sum_{\be=(v,\,u)\in\cA_c}
 \big(f(v)-e^{i\lan\t(\be),\,\vt\ran}f(u)\big), \qqq
 v\in V_c,\qq f\in\ell^2(V_c),
\]
the potential $Q$ on $\ell^2(V_c)$
has a finite rank. Here $\t(\be)\in\Z^d$ is the index of the edge $\be\in\cA_c$
defined by \er{in},\er{inf}, $V_c$ and $\cA_c$ are the vertex set and the set of oriented edges of the cylinder $\cC$, respectively; $\lan\cdot\,,\cdot\ran$ is the inner product in $\R^d$.

ii) For each $\vt\in\T^d$ the spectrum of the fiber operator
$H(\vt)$ has the form
\[\lb{strs}
\s\big(H(\vt)\big)=\s_{ac}\big(H(\vt)\big)\cup\s_{fb}\big(H(\vt)\big)
\cup\s_{p}\big(H(\vt)\big),
\]
\[\lb{strs1}
\s_{ac}\big(H(\vt)\big)=\s_{ac}\big(H_0(\vt)\big),\qqq \s_{fb}\big(H(\vt)\big)=\s_{fb}\big(H_0(\vt)\big),
\]
$\s_{p}\big(H(\vt)\big)$ is the set of all eigenvalues of
$H(\vt)$ of finite multiplicity given by
\[
\label{eq.3'}
\l_{1}(\vt)\leq\l_{2}(\vt)\leq\ldots\leq\l_{N_\vt}(\vt),\qqq
N_\vt\leq p=\rank (Q\upharpoonright V_c).
\]
\end{proposition}

\no \textbf{Remark.} The fiber Laplacian $\D(\vt)$, $\vt\in\T^d$, can
be considered as a magnetic Laplacian with a periodic
magnetic potential on the cylinder $\cC$ (see \cite{HS99a}, \cite{HS99b}, \cite{KS17}).

\

Proposition \ref{TSDI} and standard arguments (see Theorem XIII.85 in
\cite{RS78}) describe the spectrum of the guided Schr\"odinger operator
$H$. Since $H(\vt)$ is self-adjoint and analytic in $\vt\in\T^d$, each
$\l_j(\cdot)$ is a real and piecewise analytic function on the torus $\T^d$ and
creates the \emph{guided band} $\gs_j(H)$ given by
\[\lb{sgSo}
\gs_j(H)=[\l_j^-,\l_j^+]=\l_j(\T^d), \qqq j=1,\ldots,N,\qqq
N=\max\limits_{\vt\in\T^d}N_\vt\leq p.
\]
Thus, the spectrum of the guided Schr\"odinger operator $H$ on the graph $\G$ has
the form
$$
\s(H)=\bigcup_{\vt\in\T^d}\s\big(H(\vt)\big)=\s(H_0) \cup\gs(H),
$$
where $\s(H_0)$ is defined by \er{sH0} and
$$
\gs(H)=\bigcup\limits_{\vt\in\T^d}\s_{p}\big(H(\vt)\big)
=\bigcup_{j=1}^N\gs_j(H)=\gs_{ac}(H)\cup\gs_{fb}(H),
$$
$\gs_{ac}(H)$ and $\gs_{fb}(H)$ are the absolutely continuous part
and the flat band part of the guided spectrum $\gs(H)$,
respectively. An open interval between two neighboring
non-degenerate bands is called \emph{a spectral gap}. The guided
spectrum $\gs(H)$  may partly lie below the spectrum of the
unperturbed operator $H_0$, on the spectrum of $H_0$ and in the gaps
of $H_0$.

\subsection{Estimates of guided bands}
We consider the guided bands from \er{sgSo} (or their parts) below the spectrum
of the unperturbed Schr\"odinger operator $H_0$, i.e., below $\inf \s(H_0)=0$:
\[\lb{gs+}
\gs_j^o(H)=\gs_j(H)\cap(-\iy,0]\neq\varnothing,\qqq j=1,\ldots,N_g,\qqq N_g\leq N.
\]
We rewrite the sequence $Q(v),
v\in\supp(Q\upharpoonright V_c)$ in the form
\[
\lb{wtqn} 0<Q_p^\bu\le\ldots\le Q_2^\bu\le Q_1^\bu,
\]
where $Q_j^\bu=Q(v_j)$, $j=1,\ldots,p$, for some distinct vertices
$v_1,v_2,\ldots,v_p\in\supp(Q\upharpoonright V_c)$.

Proposition \ref{TSDI} and the standard perturbation theory give the estimates of the position of the bands $\gs_j^o(H)$ and their number $N_g$ by
\[\lb{esbp0}
\gs_j^o(H)\ss[-Q_j^\bu,-Q_j^\bu+\vr],\qqq N_g\geq\#\{j\in\N_p\,:\, Q_j^\bu>\vr\},
\]
where $\vr$ is defined in \er{mm} (for more details see Corollary \ref{Tloc}), $\#A$ is the number of elements of the set $A$.
In particular, this yields that the guided spectrum may lie below any fixed point and
for guided potentials $Q$ satisfying $Q_p^\bu>\vr$ and $Q_j^\bu-Q_{j+1}^\bu>\vr$
for all $j\in\N_{p-1}$ the guided spectrum of $H=H_0-Q$
consists of exactly $p$ guided bands separated by gaps.

\

In order to formulate our main result we define the set
$\cB_c=\cB/\Z^d$ of all bridges of the cylinder $\cC=(V_c,\cE_c)$ and the
modified cylinder $\cC_\gm=(V_c,\cE_c\sm\cB_c)$, which is obtained
from $\cC$ by deleting all its bridges.  We consider
the Schr\"odinger operator $h$ on the modified cylinder $\cC_\gm$:
\[
\lb{avo}
h=h_0-Q, \qqq h_0=\D_\gm+W.
\]
Here $\D_\gm$ is the Laplacian defined by \er{DOLN} on the cylinder
$\cC_\gm$ and $W$ is the restriction of the periodic potential
defined by \er{potW} to the vertex set $V_c$. The potential $Q$ has
a finite support $\{v_1,\ldots,v_p\}$ on the cylinder $\cC_\gm$.
Then the  Schr\"odinger operator $h$ has at most $p=\rank Q$ eigenvalues $\wt\m_1\leq\wt\m_2\leq\ldots$\,. Define $\m_j$ by
\[\lb{evdm}
\m_j=\min\{\wt\m_j, \inf\s_{ess}(h)\}, \qqq j=1,2,\ldots,p.
\]
We estimate the position of the guided bands $\gs_j^o(H)$ defined
by \er{gs+} in terms of the eigenvalues of the operator $h$ and
the number of bridges on the cylinder~$\cC$.

\begin{theorem}\lb{Est}
Let $H=H_0-Q$ be a guided Schr\"odinger operator. Then each guided band $\gs_j^o(H)$, $j=1,\ldots,N_g$, defined by \er{gs+} satisfies
\[\lb{esbp1}
\gs_j^o(H)\ss[\m_j,\m_j+2\b_+], \qqq \b_+=\max_{v\in V_c}\b_v,
\]
where $\b_v$ denotes the number of bridges on $\cC$ starting at the
vertex $v\in V_c$.
\end{theorem}

\no \textbf{Remarks.} 1) For most of graphs the number $\b_+=1$ and,
consequently, the length of each guided band $|\gs_j^o(H)|\leq2$,
$j=1,\ldots,N_g$. But for specific graphs $\b_+$ may be any given
integer number.

2) Note that there exists a periodic graph $\G$ such that the
inclusions \er{esbp1} become identities and $|\gs_j^o(H)|=2\b_+$, $j=1,\ldots,N_g$, for any guided potential $Q$,
see Example \ref{Ex}.

\

We discuss the guided spectrum in the case of the guided potentials large enough.

\begin{theorem}
\lb{T17} Let $H_t=\D+W-tQ$ be a guided Schr\"odinger operator, where
the coupling constant $t>0$ is large enough. Then $N_g=p$ and the
following statements hold true:

i) Let $Q^\bu_j\neq Q^\bu_k$ for some fixed $j\in\N_p$ and all
$k\in\N_p\sm \{j\}$. Then the guided band
$\gs_j(H_t)=[\l_j^-(t),\l_j^+(t)]$ satisfies
\[
\lb{Qt}
\begin{aligned}
&\l_j^\pm(t)=-tQ^\bu_j+W(v_j)+\D_j^\pm+O(1/t),\\
&|\gs_j(H_t)|=\D_j+O(1/t),
\end{aligned}
\]
as $t\to \iy$, where
\[\lb{Dpm}
\begin{aligned}
&\D_j^-=\min_{\vt\in\T^d}\D_{jj}(\vt)=\vk_{v_j}-\vk_{jj},\qqq
\D_j^+=\max_{\vt\in\T^d}\D_{jj}(\vt),\\
& \D_j=\D_j^+-\D_j^-,\qqq \b_{jj}\leq \D_j\leq2\b_{jj},
\end{aligned}
\]
for some function $\D_{jj}(\vt)$ defined by the formula \er{maDx}. Here
$\vk_{v_j}$ is the degree of the vertex $v_j\in\supp(Q\upharpoonright V_c)$, $\b_{jj}$ is the number of all oriented bridge-loops at this vertex on the cylinder $\cC$ and $\vk_{jj}$ is the number of the remaining oriented loops at $v_j$ on $\cC$.

ii) Let $Q$ be "generic", i.e., $Q_j\neq Q_k$ for all
$j,k\in\N_p$, $k\neq j$. Then the Lebesgue measure $|\gs(H_t)|$ of
the guided spectrum of the operator $H_t$ satisfies
\[
\lb{Qt1} |\gs(H_t)|=\sum\limits_{j=1}^p \D_j+O(1/t).
\]

iii) In particular, if there exists a bridge-loop at
a vertex $v_j\in\supp(Q\upharpoonright V_c)$ on the cylinder $\cC$, then $\D_j\neq0$ and the guided
band $\gs_j(H_t)$ is non-degenerate. If there are no such bridge-loops at
$v_j$, then $\D_j=0$ and $|\gs_j(H_t)|=O(1/t)$.
Moreover, if there are no bridge-loops at each vertex
$v\in\supp(Q\upharpoonright V_c)$ on $\cC$, then
$|\gs(H_t)|=O(1/t)$.
\end{theorem}

\no \textbf{Remark.} There exist periodic graphs $\G$ and guided potentials $Q$
such that the Lebesgue measure of the guided spectrum of $H=H_0-Q$ on $\G$ can be
arbitrarily large or arbitrarily small (for more details see Corollary \ref{TEg}).

\

We present the plan of our paper. In Section \ref{Sec2} we prove Proposition \ref{TSDI} about the decomposition  of guided Schr\"odinger operators into a constant fiber direct integral. In Section \ref{Sec3} we prove Theorem \ref{Est} describing the localization of the guided bands and Theorem \ref{T17} about the asymptotics of the guided bands for
large guided potentials and finally describe geometric properties of the guided spectrum for
specific graphs and guided potentials.

\section{Direct integral for guided Schr\"odinger operators}
\setcounter{equation}{0} \lb{Sec2}

\subsection{Edge indices.} In order to give a decomposition of guided Schr\"odinger
operators  into a constant fiber direct integral with a precise
representation of fiber operators we need to define {\it an edge
index}. Recall that an edge index was introduced in \cite{KS14} and
it was important to study the spectrum, effective masses of
Laplacians and Schr\"odinger operators on periodic graphs, since
fiber operators are expressed in terms of edge indices (see
\er{l2.15''}).

For any
$v\in V$ the following unique representation holds true:
\[
\lb{Dv} v=v_0+[v], \qquad v_0\in V_c,\qquad [v]\in\Z^d,
\]
where $V_c$ is the fundamental vertex set of the periodic graph
$\G$ defined by \er{fvs}. In other words, each vertex $v$ can be obtained
from a vertex $v_0\in V_c$ by the shift by a vector $[v]\in \Z^d$.
For any
oriented edge $\be=(u,v)\in\cA$ we define {\bf the edge "index"}
$\t(\be)$ as the integer vector given by
\[
\lb{in}
\t(\be)=[v]-[u]\in\Z^d,
\]
where, due to \er{Dv}, we have
$$
u=u_0+[u],\qquad v=v_0+[v], \qquad u_0,v_0\in V_c,\qquad [u],[v]\in\Z^d.
$$
We note that edges connecting vertices from the fundamental vertex set $V_c$
have zero indices.

We define a surjection $\gf_\cA:\cA\rightarrow\cA_c=\cA/\Z^d$,
which map each edge to its equivalence class. If $\be$ is an oriented edge of the graph $\G$, then there is an oriented edge
$\be_*=\gf_{\cA}(\be)$ on the cylinder $\cC=\G/\Z^d$. For the edge
$\be_*\in\cA_c$ we define the edge index $\t(\be_*)$ by
\[
\lb{inf}
\t(\be_*)=\t(\be).
\]
In other words, edge indices of the cylinder $\cC$  are
induced by edge indices of the periodic graph $\G$. Edges with nonzero indices are called \emph{bridges}.
Edge indices, generally speaking, depend on the choice of the coordinate
origin $O$ and the periods $a_1,\ldots,a_{\wt d}$ of the graph $\G$. But in a fixed coordinate system indices of the
cylinder edges are uniquely determined by \er{inf}, since
$$
\t(\be+m)=\t(\be),\qqq \forall\, (\be,m)\in\cA \ts \Z^d.
$$

\subsection{Direct integrals and an example of the guided spectrum}

We prove Proposition \ref{TSDI} about the decomposition of guided Schr\"odinger operators into a constant fiber direct integral and give an example of the guided spectrum.

\

\no \textbf{Proof of Proposition \ref{TSDI}.} i) Repeating  the
arguments from the proof of Theorem 1.1 in \cite{KS14} we obtain
\er{raz}, \er{l2.15''}, where the unitary operator
$U:\ell^2(V)\to\mH$ has the form
\[
\lb{5001} (Uf)(\vt,v)=\sum\limits_{m\in\Z^d}e^{-i\lan m,\vt\ran }
f(v+m), \qqq (\vt,v)\in \T^d\ts V_c, \qqq f\in \ell^2(V).
\]
The Hilbert space $\mH$ defined in \er{raz} is
equipped with the norm
$\|g\|^2_{\mH}=\int_{\T^d}\|g(\vt,\cdot)\|_{\ell^2(V_c)}^2\frac{d\vt
}{(2\pi)^d}$\,, where the function $g(\vt,\cdot)\in\ell^2(V_c)$ for almost all
$\vt\in\T^d$. The fiber potential $Q$ is a finite rank operator,
since $\supp(Q\upharpoonright V_c)=\{v_1,\ldots,v_p\}\ss V_c$ is
finite.

ii) For each $\vt\in\T^d$ the unperturbed fiber operator $H_0(\vt)$ is $\Z^{\wt
d-d}$-periodic. Then, using standard arguments (see
Theorem XIII.85 in \cite{RS78}), we obtain that the spectrum of $H_0(\vt)$ has the form
$$
\s\big(H_0(\vt)\big)=\s_{ac}\big(H_0(\vt)\big)\cup
\s_{fb}\big(H_0(\vt)\big).
$$
Since the fiber potential $Q$ is an operator of
rank $p=\#\supp(Q\upharpoonright V_c)$, for each
$\l\in\s_{fb}\big(H_0(\vt)\big)$ there exists
a corresponding eigenfunction with a finite support (see, e.g., Theorem 4.5.2 in \cite{BK13}) not intersecting with
$\supp(Q\upharpoonright V_c)$. Thus, $\l\in\s_{fb}\big(H(\vt)\big)$ and vice versa.
Then for each $\vt\in\T^d$
the spectrum $\s\big(H(\vt)\big)$ of the fiber operator
$H(\vt)$ is given by \er{strs}, where $\s_{ac}\big(H(\vt)\big)$,
$\s_{fb}\big(H(\vt)\big)$ satisfy \er{strs1} and
$\s_{p}\big(H(\vt)\big)$ consists of $N_\vt\leq p$ eigenvalues \er{eq.3'}.
\qq $\BBox$

\setlength{\unitlength}{1.0mm}
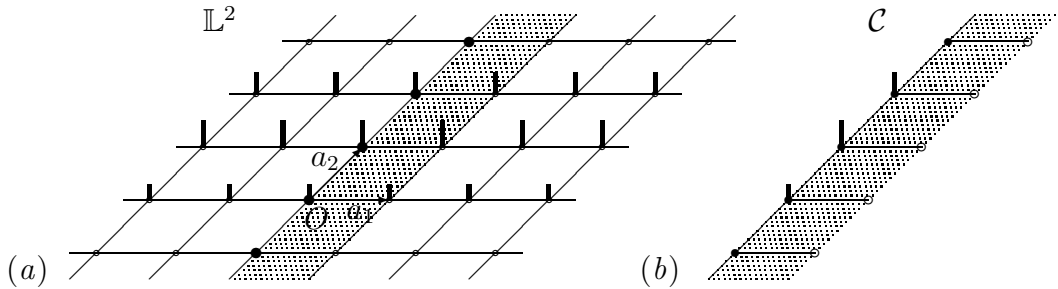
\begin{figure}[h]
\centering

\unitlength 0.7mm 
\linethickness{0.4pt}
\ifx\plotpoint\undefined\newsavebox{\plotpoint}\fi 
\begin{picture}(200,60)(0,0)

\put(-12,5){(\emph{a})} \put(25,52){$\dL^2$}

\put(0,10){\line(1,0){85.00}} \put(10,20){\line(1,0){85.00}}
\put(20,30){\line(1,0){85.00}} \put(30,40){\line(1,0){85.00}}
\put(40,50){\line(1,0){85.00}}

\put(0,5){\line(1,1){50.00}} \put(15,5){\line(1,1){50.00}}
\put(30,5){\line(1,1){50.00}} \put(45,5){\line(1,1){50.00}}
\put(60,5){\line(1,1){50.00}} \put(75,5){\line(1,1){50.00}}

\put(5,10){\circle{1}} \put(20,10){\circle{1}}
\put(35,10){\circle*{2}} \put(50,10){\circle{1}}
\put(65,10){\circle{1}} \put(80,10){\circle{1}}

\put(15,20){\circle{1}} \put(30,20){\circle{1}}
\put(45,20){\circle*{2}} \put(45,20){\vector(1,0){15.00}}
\put(45,20){\vector(1,1){10.00}} \put(44,14.5){$O$}
\put(52,16.5){$a_1$} \put(45.3,26.5){$a_2$} \put(60,20){\circle{1}}
\put(75,20){\circle{1}} \put(90,20){\circle{1}}

\bezier{60}(31,5)(56,30)(81,55) \bezier{60}(32,5)(57,30)(82,55)
\bezier{60}(33,5)(58,30)(83,55) \bezier{60}(34,5)(59,30)(84,55)
\bezier{60}(35,5)(60,30)(85,55) \bezier{60}(36,5)(61,30)(86,55)
\bezier{60}(37,5)(62,30)(87,55) \bezier{60}(38,5)(63,30)(88,55)
\bezier{60}(39,5)(64,30)(89,55) \bezier{60}(40,5)(65,30)(90,55)
\bezier{60}(41,5)(66,30)(91,55) \bezier{60}(42,5)(67,30)(92,55)
\bezier{60}(43,5)(68,30)(93,55) \bezier{60}(44,5)(69,30)(94,55)
\bezier{60}(45,5)(70,30)(95,55)

\linethickness{1.5pt} \put(15,20){\line(0,1){3.00}}
\put(30,20){\line(0,1){3.00}} \put(45,20){\line(0,1){3.00}}
\put(60,20){\line(0,1){3.00}} \put(75,20){\line(0,1){3.00}}
\put(90,20){\line(0,1){3.00}}

\put(25,30){\line(0,1){5.00}} \put(40,30){\line(0,1){5.00}}
\put(55,30){\line(0,1){5.00}} \put(70,30){\line(0,1){5.00}}
\put(85,30){\line(0,1){5.00}} \put(100,30){\line(0,1){5.00}}

\put(35,40){\line(0,1){4.00}} \put(50,40){\line(0,1){4.00}}
\put(65,40){\line(0,1){4.00}} \put(80,40){\line(0,1){4.00}}
\put(95,40){\line(0,1){4.00}} \put(110,40){\line(0,1){4.00}}
\linethickness{0.4pt}

\put(25,30){\circle{1}} \put(40,30){\circle{1}}
\put(55,30){\circle*{2}} \put(70,30){\circle{1}}
\put(85,30){\circle{1}} \put(100,30){\circle{1}}

\put(35,40){\circle{1}} \put(50,40){\circle{1}}
\put(65,40){\circle*{2}} \put(80,40){\circle{1}}
\put(95,40){\circle{1}} \put(110,40){\circle{1}}

\put(45,50){\circle{1}} \put(60,50){\circle{1}}
\put(75,50){\circle*{2}} \put(90,50){\circle{1}}
\put(105,50){\circle{1}} \put(120,50){\circle{1}}

\put(150,52){$\cC$} \put(125,10){\line(1,0){15.00}}
\put(135,20){\line(1,0){15.00}} \put(145,30){\line(1,0){15.00}}
\put(155,40){\line(1,0){15.00}} \put(165,50){\line(1,0){15.00}}

\put(120,5){\line(1,1){50.0}} \bezier{60}(135,5)(160,30)(185,55)

\bezier{60}(134,5)(159,30)(184,55)
\bezier{60}(133,5)(158,30)(183,55)
\bezier{60}(132,5)(157,30)(182,55)
\bezier{60}(131,5)(156,30)(181,55)
\bezier{60}(130,5)(155,30)(180,55)
\bezier{60}(129,5)(154,30)(179,55)
\bezier{60}(128,5)(153,30)(178,55)
\bezier{60}(127,5)(152,30)(177,55)
\bezier{60}(126,5)(151,30)(176,55)
\bezier{60}(125,5)(150,30)(175,55)
\bezier{60}(124,5)(149,30)(174,55)
\bezier{60}(123,5)(148,30)(173,55)
\bezier{60}(122,5)(147,30)(172,55)
\bezier{60}(121,5)(146,30)(171,55)
\bezier{60}(120,5)(145,30)(170,55)

\linethickness{1.5pt} \put(135,20){\line(0,1){3.00}}
\put(145,30){\line(0,1){5.00}} \put(155,40){\line(0,1){4.00}}
\linethickness{0.4pt}

\put(125,10){\circle*{1.5}} \put(140,10){\circle{1.5}}

\put(135,20){\circle*{1.5}} \put(150,20){\circle{1.5}}

\put(145,30){\circle*{1.5}} \put(160,30){\circle{1.5}}

\put(155,40){\circle*{1.5}} \put(170,40){\circle{1.5}}

\put(165,50){\circle*{1.5}} \put(180,50){\circle{1.5}}

\put(107,5){(\emph{b})}
\end{picture}

\vspace{-0.3cm} \caption{\footnotesize  \emph{a}) The square lattice
$\dL^2$; the vertices from the fundamental vertex set $V_c$ are big
black points; the strip $\cS$ is shaded;\; \emph{b}) the cylinder
$\cC=\dL^2/\Z$ (the edges of the strip are identified). The values
of the guided potential $Q$ are along the vertical axis.}
\label{SLp}
\end{figure}
We consider a simple example of the guided spectrum.
\begin{example}\lb{Ex}
Let $\dL^2=(V,\cE)$ be the square lattice, where the vertex set
$V=\Z^2$ and the edge set $\cE=\big\{(m,m+a_s), \;
\forall\,m\in\Z^2, s=1,2\big\}$ and $a_1=(1,0)$, $a_2=(0,1)$, see
Fig.\ref{SLp}.a.
Our Laplacian $\D$ defined by \er{DOLN} has
 the form
\[\lb{Lao}
(\D f)(m)=4f(m)-\sum_{|m-k|=1}f(k), \qqq f\in\ell^2(\Z^2), \qqq
m\in\Z^2.
\]
We consider a Schr\"odinger operator  $H=\D-Q$ on $\Z^2$, where
$Q\ge 0$ is a guided potential such that
$$ (Q f)(m)=Q(m_2)f(m),
\qqq m=(m_1,m_2),
$$
and $Q(m_2)$ is finitely supported on $\Z$.
\end{example}
Due to Proposition
\ref{TSDI}, the guided operator  $H=\D-Q$ on the lattice $\Z^2$ has
the decomposition \er{raz} into a constant fiber direct integral,
where the fiber Schr\"odinger operator $H(\vt)=\D(\vt)-Q$ acts on
$f\in\ell^2(\Z)$ and is given by
\[
\label{l2.13}
\begin{aligned}
& H(\vt)=2(1-\cos\vt)+h, \\
&  (hf)(n)= 2f(n)-f(n+1)-f(n-1)-Q(n)f(n), \qqq n\in\Z,
\end{aligned}
\]
for all $\vt\in \T=(-\pi,\pi]$. The operator $h$ is well studied,
see \cite{T89} for a large class of $Q$, and see \cite{K11} in the
case of finitely supported $Q$, where the inverse problem is solved.
It is well known that  the spectrum of the operator $h$ consists of
an absolutely continuous part $[0,4]$ plus a finite number of simple
eigenvalues
$$
\m_1<\m_2<\ldots<\m_N<0, \qqq N\leq p, \qqq p=\rank Q.
$$
Then, each fiber operator $H(\vt)$, $\vt\in\T$, has at most $p$
simple eigenvalues:
$$
\l_1(\vt)<\l_2(\vt)<\ldots<\l_N(\vt)<2-2\cos\vt, \qq {\rm where }\qq
\l_j(\vt)=\m_j+2-2\cos\vt,\qq j\in\N_N.
$$
Thus, the definition \er{sgSo} of the guided bands $\gs_j(H)$ yields
$$
\gs_j(H)=\l_j(\T)=[\m_j,\m_j+4],\qqq |\gs_j(H)|=4, \qqq j\in\N_N,
\qqq N\leq p.
$$
On the other hand, the number $\b_+$ defined in \er{esbp1} is equal to 2, since
for each vertex of the cylinder $\cC$ there are two bridges starting at this vertex.
Thus, on the graph $\dL^2$ the inclusions \er{esbp1} become identities. \qq
\BBox

\section{Proof of the main results}
\setcounter{equation}{0}
\lb{Sec3}

In this section we
prove Theorem \ref{Est} about
the position of the guided bands and Theorem \ref{T17} about the asymptotics of the guided bands for large guided potentials. We also describe geometric properties of the guided spectrum for
specific graphs and guided potentials (see Corollary \ref{TEg}).

\subsection{Estimates for the guided spectrum.} Denote by $m_\pm(\vt)$
the upper and lower endpoints of the spectrum of the unperturbed
fiber operator $H_0(\vt)$:
\[\lb{mpmm}
m_-(\vt)=
\inf\s\big(H_0(\vt)\big),\qqq m_+(\vt)= \sup\s\big(H_0(\vt)\big).
\]
Then \er{mm} yields
\[
\lb{mvt}
\min_{\vt\in\T^d}m_-(\vt)=0,\qqq \max_{\vt\in\T^d}m_+(\vt)=\vr.
\]

We need a simple estimate for eigenvalues of bounded self-adjoint operators \cite{RS78}:
\emph{Let $A,B$ be bounded self-adjoint
operators in a Hilbert space $\cH$ and let
$\l_j(A)=\min\big\{\wt\l_j(A),\inf\s_{ess}(A)\big\}$,
$j=1,2,\ldots$\,, where $\wt \l_1(A)\leq\wt \l_2(A)\leq\ldots$ are
the eigenvalues of $A$. Then}
\[\lb{wnf1}
\l_j(A)+\inf\s(B)\leq\l_j(A+B)\leq\l_j(A)+\sup\s(B), \qqq j=1,2,3,\ldots\,.
\]
The following simple corollary about the position of the guided bands $\gs_j^o(H)$  defined by \er{gs+} is a direct consequence of Proposition \ref{TSDI}.

\begin{corollary}\lb{Tloc}
Let $H=H_0-Q$ be a guided Schr\"odinger operator and let $\vr$ be
defined in \er{mm}. Then each guided band $\gs_j^o(H)$,
$j=1,\ldots,N_g$, and their number $N_g$ satisfy
\[\lb{esbp}
\gs_j^o(H)\ss[-Q_j^\bu,-Q_j^\bu+\vr],
\]
\[\lb{esNg}
N_g\geq\#\{j\in\N_p\,:\, Q_j^\bu>\vr\}.
\]
\end{corollary}

\no \textbf{Proof.} The fiber Schr\"odinger operator $H(\vt)$
is given by $H(\vt)=H_0(\vt)-Q$, $Q\geq0$. Then, due to \er{wnf1}, \er{mpmm} and \er{mvt},
for each $\vt\in\T^d$ the eigenvalues $\l_j(\vt)$ of $H(\vt)$ below its essential spectrum satisfy
\[\lb{com1}
-Q_j^\bu\leq m_-(\vt)-Q_j^\bu\leq\l_j(\vt)\leq m_+(\vt)-Q_j^\bu\leq
\vr-Q_j^\bu,
\]
which yields \er{esbp}. Let $Q_j^\bu>\vr$ for some $j=1,\ldots,p$.
Then, due to \er{com1}, $\l_j(\vt)<0$ for all $\vt\in\T^d$. Thus,
$\l_j$ creates the guided band $\gs_j^o(H)=\l_j(\T^d)$. This yields
\er{esNg}. \qq \BBox

\

\no \textbf{Proof of Theorem \ref{Est}.} We rewrite the fiber operator $H(\vt)$,
$\vt\in\T^d$, defined by \er{raz}, \er{l2.15''} in the form:
\[
\label{eq.1}
H(\vt)=h+\D_\b(\vt),\qqq h=\D_\gm+W-Q,
\]
\[\label{Dbt}
\big(\D_\b(\vt)f\big)(v)=\sum_{\be=(v,\,u)\in\cA_c \atop \t(\be)\neq0}
\big(f(v)-e^{i\lan\t(\be),\,\vt\ran}f(u)\big), \qqq
 v\in V_c,
\]
where $\t(\be)\in\Z^d$ is the index of the edge $\be\in\cA_c$
defined by \er{in}, \er{inf}. Each operator $\D_\b(\vt)$, $\vt\in\T^d$,
is the magnetic Laplacian on the graph $\cC_\b=(V_c,\cB_c)$ and the
degree of each vertex $v\in V_c$ on $\cC_\b$ is equal to the number
$\b_v$ of all bridges starting at $v$. Then the  spectrum of
$\D_\b(\vt)$ satisfies
$$
\s\big(\D_\b(\vt)\big)\ss[0,2\b_+], \qqq \forall\, \vt\in\T^d
$$
(see, e.g., \cite{HS99a}, \cite{HS99b}) and, due to \er{wnf1},  we have
$\m_j\le\l_j(\vt)\le\m_j+2\b_+$ for each $\vt\in\T^d$ and each $j=1,\ldots,N_g$, which yields \er{esbp1}.
\qq \BBox


\subsection {Proof of Theorem \ref{T17}.} We consider the guided Schr\"odinger operator $H_t=H_0-tQ$, $Q\geq0$. If $t>0$ is large enough, then $tQ_j^\bu>\vr$ for each
$j=1,\ldots,p$, where $\vr$ is defined in \er{mm}. Then, due to \er{esNg}
and the inequality $N_g\leq p$, we have $N_g=p$.

i) We rewrite the fiber operator $H_t(\vt)$, $\vt\in\T^d$, for the
guided operator $H_t$ in the form
$$
H_t(\vt)=H_0(\vt)-tQ=t K_t(\vt),\qqq K_t(\vt)=-Q+\ve H_0(\vt),\qqq
\ve=\frac1t\,.
$$
We denote the eigenvalues of the operator $K_t(\vt)$ below its
essential spectrum by
\[
\label{ebesH} E_1(\vt,t)\leq E_2(\vt,t)\leq\ldots\leq E_p(\vt,t),
\qqq \vt\in\T^d.
\]
 For each vertex $u\in\supp(Q\upharpoonright V_c)$ we define the
function $f_u\in\ell^2(V_c)$ by $f_u(v)=\d_{uv}$, where $\d_{uv}$ is the Kronecker delta. Then the
eigenvalue $E_j(\vt,t)$ of the operator $K_t(\vt)$  has the
following asymptotics:
\[\lb{ass0}
E_j(\vt,t)=-Q^\bu_j+\ve\,\lan f_{v_j},H_0(\vt)f_{v_j}\ran_{V_c}+\ve^2\sum_{k=1 \atop k\neq
j}^p\frac{|\lan f_{v_k},H_0(\vt)f_{v_j}\ran_{V_c}|^2}{Q^\bu_k-Q^\bu_j}+O(\ve^3)
\]
(see pp. 7--8 in \cite{RS78}) uniformly in $\vt\in\T^d$ as
$t\rightarrow\infty$, where $\lan\cdot\,,\cdot\ran_{V_c}$ denotes the inner product in $\ell^2(V_c)$.
Using the identity $H_0(\vt)=\D(\vt)+W$ and the formula
\er{l2.15''} for the fiber Laplacian $\D(\vt)$, we rewrite the asymptotics \er{ass0} in the form
\[
\lb{ass1}
E_j(\vt,t)=-Q^\bu_j+\ve\big(\D_{jj}(\vt)+W(v_j)\big)+O(\ve^2),
\]
where $\D_{jj}(\vt)$ is defined by
\[
\lb{maDx} \D_{jj}(\vt)=\vk_{v_j}-\sum_{\be=(v_j,\,v_j)\in\cA_c}
\cos\lan\t(\be),\vt\ran.
\]
Here $\vk_v$ is the degree of the
vertex $v$,  the vertices
$v_1,v_2,\ldots,v_p\in\supp(Q\upharpoonright V_c)$ and
$\t(\be)\in\Z^d$ is the edge index defined by \er{in}, \er{inf}.
This yields the asymptotics of
the eigenvalue $\l_j(\vt,t)$ of the operator $H_t(\vt)$:
\[\lb{ass}
\l_j(\vt,t)=t\,E_j(\vt,t)=-tQ^\bu_j+\D_{jj}(\vt)+W(v_j)+ O(1/t).
\]
Using this asymptotics for $\l_j^-(t)=\min\limits_{\vt\in\T^d}\l_j(\vt,t)$ and $\l_j^+(t)=\max\limits_{\vt\in\T^d}\l_j(\vt,t)$, we
obtain
$$
\begin{aligned}
&\l_j^\pm(t)=-tQ^\bu_j+\D_j^\pm +W(v_j)+O(1/t),\\
\end{aligned}
$$
where $\D_j^\pm$ are defined in \er{Dpm}. Since
$\gs_j(H_t)=[\l_j^-(t),\l_j^+(t)]$, the asymptotics \er{ass} also
gives the second formula in \er{Qt}. Using \er{maDx} we rewrite the constant $\D_j$ defined in \er{Dpm} in the form
\[\lb{esCj}
\D_j=\max_{\vt\in\T^d}\Omega_j(\vt)-\min_{\vt\in\T^d}\Omega_j(\vt)=
\b_{jj}-\min_{\vt\in\T^d}\Omega_j(\vt),\qq \Omega_j(\vt)=\sum_{\be=(v_j,\,v_j)\in\cA_c \atop \t(\be)\neq0}
\cos\lan\t(\be),\vt\ran.
\]
Using the identity $\int\limits_{\T^d}\cos\lan\t(\be),\vt\ran\,d\vt=0$ for each $\t(\be)\neq0$, we obtain $-\b_{jj}\leq\min\limits_{\vt\in\T^d}\Omega_j(\vt)\leq0$ and \er{esCj} yields $\b_{jj}\leq \D_j\leq 2\b_{jj}$.

ii) If the guided potential $Q$ is generic, then summing the second asymptotics in \er{Qt} over $j=1,\ldots,p$ we obtain \er{Qt1}.

iii) Let on the cylinder $\cC$ there
exist a bridge-loop at a vertex $v_j\in\supp(Q\upharpoonright V_c)$, i.e., a loop $\be$ with
$\t(\be)\neq0$. This yields that the function $\D_{jj}$ defined by
\er{maDx} is not constant, i.e., $\D_j\neq0$, and, due to the second asymptotics in
\er{Qt}, the guided band $\gs_j(H_t)$ is non-degenerate. If there
are no bridge-loops at $v_j$ on $\cC$, then $\D_{jj}$ is constant, i.e., $\D_j=0$, and, using the second
asymptotics in \er{Qt}, we obtain $|\gs_j(H_t)|=O(1/t)$.

If for each $v\in\supp(Q\upharpoonright V_c)$ there are no bridge-loops
at $v$ on $\cC$, then
$\D_{j}=0$ for each $j\in \N_p$ and the asymptotics \er{Qt1}
takes the form $|\gs(H_t)|=O(1/t)$. \qq $\BBox$

\

\no \textbf{Remarks.} i) The set $\cA_c$ of all oriented edges  of the
cylinder $\cC$ is infinite, but the sum in \er{maDx} is taken over
a finite (maybe empty) set of edge-loops at the vertex $v_j$.

ii) If all bridge-loops at some vertex $v_j\in\supp(Q\upharpoonright V_c)$ on
$\cC$ have linearly independent indices, then $\D_j$ defined in
\er{Dpm} is equal to $2\b_{jj}$.

\

Now we describe geometric properties of the guided spectrum for
specific graphs and guided potentials.

\begin{corollary}\lb{TEg}
Let  $H=H_0-Q$ be a guided Schr\"odinger operator. Then the
following statements hold true.

i) Let there exist a bridge-loop on the cylinder $\cC$.
Then for any constant $C>0$ there exists a guided potential $Q\geq0$
such that the Lebesgue measure of the guided spectrum
$\gs(H)$ satisfies $|\gs(H)|>C$ and all guided
bands are non-degenerate.

ii) Let there exist a vertex $v\in V_c$ such that there is no
bridge-loop at $v$ on $\cC$. Then for
any small $\ve>0$ there exists a guided potential $Q\geq0$ such that
the Lebesgue measure of the guided spectrum $\gs(H)$
satisfies $|\gs(H)|<\ve$.

\end{corollary}

\no \textbf{Remark.} An example of a cylinder $\cC$ with
bridge-loops is shown in Fig.\ref{SLp}.\emph{b}.

\

\no \textbf{Proof.} We consider guided Schr\"odinger operators $H_t=H_0-tQ$, $Q\geq0$ with the coupling constant $t>0$ large enough.
Due to Theorem \ref{T17}, the number $N_g$ of the guided bands
$\gs_j^o(H_t)=\gs_j(H_t)$ is equal to~$p$.

i) Let $\be_{11},\ldots,\be_{1q}\in\cA_c$ be all bridge-loops
at some vertex $v_1\in V_c$ and $\t(\be_{11}),\ldots,\t(\be_{1q})$
be their indices defined by \er{in}, \er{inf}. Due to the periodicity of the cylinder $\cC$,
at each vertex $v_j=v_1+(j-1)a_{\wt d}$\,, $j\in\Z$, there exist exactly
$q$ bridge-loops $\be_{j1},\ldots,\be_{jq}\in\cA_c$
with the same indices $\t(\be_{j1}),\ldots,\t(\be_{jq})$. Then we have
\[\lb{D1p}
\D_{11}(\vt)=\D_{22}(\vt)=\ldots,
\]
where $\D_{jj}$ is defined by \er{maDx}.

Let $Q$ be a "generic" guided potential, i.e., $Q_j\neq Q_k$ for all $j,k\in\N_p$, $j\neq k$,  with $\supp(Q\upharpoonright
V_c)=\{v_1,v_2,\ldots,v_p\}$ for some $p\in\N$.
Then, due to Theorem \ref{T17}, for $t$ large enough the guided bands $\gs_j(H_t)$, $j=1,\ldots,p$, satisfy
$$
|\gs_j(H_t)|=\D_j+O(1/t),
$$
where $\D_j$ is defined in \er{Dpm} and, due to \er{D1p},
the constant $\D_j\neq0$ is the same for all $j=1,\ldots,p$.
This yields that all guided bands are non-degenerate and
\[\lb{leme}
|\gs(H_t)|=\sum\limits_{j=1}^p\big|\gs_j(H_t)\big|=p\,\D_1+O(1/t).
\]
Choosing $p>\dfrac C{\D_1}$\,, we obtain $|\gs(H_t)|>C$ for $t$ large
enough.

ii) Let there exist a vertex $v\in V_c$ such that there is no bridge-loop
at $v$ on $\cC$. We consider the
guided potential $Q$ with $\supp(Q\upharpoonright V_c)=\{v\}$. Due
to Theorem \ref{T17}, for $t$ large enough the guided spectrum of
$H_t$ consists of exactly one guided band $\gs_1(H_t)$ and the length
of this guided band satisfies $|\gs_1(H_t)|=O(1/t)$. Thus, for any
small $\ve>0$ there exists $t>0$ such that $|\gs(H_t)|<\ve$. \qq
\BBox


\medskip

\footnotesize
 \textbf{Acknowledgments. \lb{Sec8}}  Our study was
supported by the RSF grant  No. 15-11-30007.

\end{document}